\documentclass[11pt]{amsart}
\usepackage{hyperref}
\usepackage{amsmath,array,,color}
\usepackage{amssymb,tabularx}
\usepackage{listings}
\usepackage{fancyvrb}
\usepackage{multirow}
\usepackage{cancel}
\usepackage{fancyvrb}
\usepackage{multirow}
\usepackage{cancel}
\usepackage{listings}
\usepackage{tikz}
\usetikzlibrary{positioning,arrows,shapes,decorations.pathmorphing,shadows}
\usepackage{tocvsec2}
\usepackage{booktabs}
\usepackage{float}

\tikzstyle{every picture}+=[remember picture]

\tikzstyle{colore1}=[ball color=SandyBrown] 
\tikzstyle{colore2}=[ball color=White] 
\tikzstyle{coloreno}=[ball color=Red] 
\tikzstyle{colore3}=[ball color=LightGrey] 
\tikzstyle{colore4}=[ball color=LightGrey] 
\tikzstyle{freccia}=[thick] 

\usepackage[a4paper]{geometry}

\newcommand{\shortform}[1]{\ensuremath{{\scriptstyle{#1}}}}
\newcommand{\boldshortform}[1]{\ensuremath{{\scriptstyle\boldsymbol{#1}}}}

\newcommand{\Id}{\mathop{\mathrm{Id}}}

\newcommand{\Spin}[1]{\ensuremath{\text{\upshape\rmfamily Spin}(#1)}}

\newcommand{\SO}[1]{\ensuremath{\text{\upshape\rmfamily SO}(#1)}}
\newcommand{\U}[1]{\ensuremath{\text{\upshape\rmfamily U}(#1)}}
\newcommand{\Gtwo}{\ensuremath{\text{\upshape\rmfamily G}_2}}
\newcommand{\SU}[1]{\ensuremath{\text{\upshape\rmfamily SU}(#1)}}

\newcommand{\Sp}[1]{\ensuremath{\text{\upshape\rmfamily Sp}(#1)}}

\newcommand{\ug}{\;\shortstack{{\tiny\upshape def}\\=}\;}
\newcommand{\snform}{\ensuremath{\Phi}}
\newcommand{\spinform}[1]{\ensuremath{\Phi_{\Spin{#1}}}}

\renewcommand{\Im}{\mathop{\mathrm{Im}}}

\newcommand{\form}[1]{\ensuremath{\snform_{#1}}}

\newcommand{\CC}{\mathbb{C}}   
\newcommand{\HH}{\mathbb{H}}   
\newcommand{\RR}{\mathbb{R}}

\newcommand{\OO}{\mathbb{O}}

\newcommand{\CP}[1]{\mathbb{C}P^{#1}}
\newcommand{\OP}[1]{\mathbb{O}P^{#1}}

\newcommand{\HP}[1]{\mathbb{H}P^{#1}}

\newcommand{\I}{\mathcal{I}} 
\newcommand{\J}{J} 

\numberwithin{equation}{section}

\newtheorem{te}{Theorem}[section]
\newtheorem*{te*}{Theorem}

\theoremstyle{definition}
\newtheorem{de}[te]{Definition}

\theoremstyle{remark}
\newtheorem{re}[te]{Remark}

\begin{document}

\title[Clifford systems, Clifford structures]{Clifford systems, Clifford structures, \\
and their canonical differential forms}

\subjclass[2010]{Primary 53C26, 53C27, 53C38.}
\keywords{Octonions, Clifford system, Clifford structure, calibration, canonical form.}

\author[K. B. Boydon]{Kai Brynne M. Boydon}
\address{Institute of Mathematics\\ University of the Philippines Diliman, Philippines}
\email{\href{mailto:kbboydon@math.upd.edu.ph}{kbboydon@math.upd.edu.ph}}

\author[P. Piccinni]{Paolo Piccinni}
\address{Dipartimento di Matematica\\ Sapienza Universit\`a di Roma\\
Piazzale Aldo Moro 2, I-00185, Roma, Italy 
}
\email{\href{mailto:piccinni@mat.uniroma1.it}{piccinni@mat.uniroma1.it}}

\thanks{The first author was supported by University of the Philippines OVPAA Doctoral Fellowship. Part of the present work was done during her visit at Sapienza Universit\`a di Roma in the academic year 2018-19, and she thanks Sapienza University and Department of Mathematics "Guido Castelnuovo" for hospitality.}

\thanks{The second author was supported by the group GNSAGA of INdAM, by the PRIN Project of MIUR ``Real and Complex Manifolds: Topology, Geometry and holomorphic dynamics'' and by Sapienza Universit\`a di Roma Project ``Polynomial identities and combinatorial methods in algebraic and geometric structures''}


\begin{abstract}
A comparison among different constructions in $\HH^2 \cong \RR^8$ of the quaternionic $4$-form $\Phi_{\Sp{2}\Sp{1}}$ and of the Cayley calibration $\Phi_{\Spin{7}}$ shows that one can start for them from the same collections of "K\"ahler 2-forms", entering both in quaternion K\"ahler and in $\Spin{7}$ geometry. This comparison relates with the  notions of even Clifford structure and of Clifford system. Going to dimension $16$, similar constructions allow to write explicit formulas in $\RR^{16}$ for the canonical $4$-forms $\form{\Spin{8}}$ and  $\form{\Spin{7}\U{1}}$, associated with Clifford systems related with the subgroups $\Spin{8}$ and  $\Spin{7}\U{1}$ of $\SO{16}$. We characterize the calibrated $4$-planes of the $4$-forms $\form{\Spin{8}}$ and  $\form{\Spin{7}\U{1}}$, extending in two different ways the notion of Cayley $4$-plane to dimension $16$.
\end{abstract}

\maketitle


\section{Introduction}\label{sec:intro}

In 1989 R. Bryant and R. Harvey defined the following calibration, of interest in hyperk\"ahler geometry \cite{br-ha}:
\[
\Phi_K= - \frac{1}{2} \omega_{R_i}^2 - \frac{1}{2} \omega_{R_j}^2+\frac{1}{2}\omega_{R_k}^2 \in \Lambda^4 \HH^n.
\] 
In this definition, $(\omega_{R_i},\omega_{R_j},\omega_{R_k})$ are the K\"ahler 2-forms of the hypercomplex structure $(R_i,R_j,R_k)$, defined by multiplications on the right by unit quaternions $(i,j,k)$ on the space $\RR^{4n} \cong \HH^n$. 

When $n=2$, the Bryant-Harvey calibration $\Phi_K$ relates with $\Spin{7}$ geometry. This is easily recognized by using the map
\[
L: \HH^2 \rightarrow \OO, \qquad L(h_1,h_2)=h_1 +(kh_2 \bar k)e\in\OO,
\] 
from pairs of quaternions to octonions, that yields the identity
\begin{equation}
\label{br-ha}
L^*\spinform{7} = \Phi_K.
\end{equation}
Here $\spinform{7} \in \Lambda^4 \RR^8$ is the $\Spin{7}$ 4-form, or \emph{Cayley calibration}, studied since the R. Harvey and H. B. Lawson's foundational paper \cite{ha-la}, and defined through the scalar product and the double cross product of $\RR^8 \cong \OO$:
\[
\spinform{7}(x,y,z,w) \; = \;  <x \,, \,y \times z \times w> \; = \;  <x\,,\, y (\bar z w)>,
\]
assuming here orthogonal $y,z,w \in \OO$.

The present paper collects some of the results in the first author Ph.D. thesis \cite{Boy}, inspired from viewing formula \eqref{br-ha} as a way of constructing the Cayley calibration $\spinform{7}$ through the 2-forms $\omega_{R_i},\omega_{R_j},\omega_{R_k}$. As well known, by summing the squares of the latter 2-forms one gets another remarkable calibration, namely the quaternionic right 4-form $\Omega_R$. Thus $\omega_{R_i},\omega_{R_j},\omega_{R_k}$, somehow building blocks for quaternionic geometry, enter also in $\Spin{7}$ geometry.

A first result is the following Theorem \ref{Thm 1.1}, a kind of "other way around" of formula \eqref{br-ha}. To state it, recall that the Cayley calibration $\spinform{7}$ can also be constructed as sum of squares of "K\"ahler 2-forms" associated with complex structures on $\RR^8 \cong \OO$, defined by the unit octonions. In fact, cf. \cite[Prop.10]{pp1}: 
\begin{equation}\label{spin7}
\spinform{7}\, =\, -\frac{1}{6}(\phi_{i}^2+\phi_{j}^2+\dots+\phi_{h}^2)\, =\, \frac{1}{6}(\varphi_{ij}^2+\phi_{ik}^2+\varphi_{ik}^2+\dots+\varphi_{gh}^2).
\end{equation}
Here $\phi_{i}, \phi_{j}, \dots ,\phi_{h}$ are the K\"ahler 2-forms associated with the 7 complex structures $R_i, R_j, \dots , R_h$ on $\RR^8 \cong \OO$, the right multiplications by the unit octonions $i,j,k,e,f,g,h$, and $\phi_{ij}, \phi_{ik}, \dots,\phi_{gh}$ are the K\"ahler 2-forms associated with the 21 complex structures $R_{ij}=R_i \circ R_j, R_{ik} =R_i \circ R_k, \dots , R_{gh}= R_g \circ R_h$, compositions of them. 

\begin{te}\label{Thm 1.1} The right quaternionic 4-form $\Omega_R \in \Lambda^4 \HH^2$ can be obtained from the the K\"ahler forms $\phi_i, \phi_j, \dots ,\phi_h$ associated with the complex structures $R_i, R_j, \dots ,R_h$ as:
  \[
 \Omega_R= 2[  \phi_i^2+  \phi_j^2+ \phi_k^2 - \phi_e^2 - \phi_f^2-\phi_g^2-  \phi_h^2].
  \]
Moreover, by selecting any five out of the seven $(J_1=R_i,J_2=R_j, \dots,J_7=R_h)$ and by looking at the matrix $\zeta=(\zeta_{\alpha \beta}) \in \frak{so}(5)$ of K\"ahler 2-forms of their compositions $J_{\alpha \beta} =  J_\alpha \circ J_\beta$, one can get the left quaternionic 4-form $\Omega_L$ as
  \[
 \Omega_L= -\frac{1}{2} \sum_{\alpha < \beta} \zeta_{\alpha \beta}^2,
  \]
up to a permutation or change of signs of some coordinates in $\RR^8$.
\end{te}

On the same direction as in Bryant-Harvey's formula \eqref{br-ha}, a similar result is the following (cf. Section \ref{sec:cayley2} for more details):

\begin{te}\label{Thm 1.2} The Cayley calibration $\spinform{7} \in \Lambda^4 \RR^8$ can be obtained from the K\"ahler 2-forms $\eta_{\alpha\beta}$ $(1\leq \alpha < \beta \leq 5)$ associated to complex structures $J^L_{\alpha\beta} = \I^L_\alpha \circ \I^L_\beta$, where $\I^L_1, \dots \I^L_5$ are anti-commuting self-dual involutions in $\RR^8$. Namely:
  \[
\Phi_{\mathrm{Spin(7)}}= \frac{1}{4} \big[ \eta_{12}^2+\eta_{13}^2+\eta_{24}^2+\eta_{34}^2-\eta_{23}^2-\eta_{14}^2-\eta_{15}^2-\eta_{25}^2-\eta_{35}^2-\eta_{45}^2 \big],
  \]
and on the other hand one can get the right quaternion K\"ahler $\Omega_R$ as: 
\[
\Omega_R\; =\;  -\frac{1}{2} \big[ \eta_{12}^2+\eta_{13}^2+\eta_{24}^2+\eta_{34}^2+\eta_{23}^2+\eta_{14}^2+\eta_{15}^2+\eta_{25}^2+\eta_{35}^2+\eta_{45}^2 \big].
\]
\end{te}
\bigskip

Moving to dimension 16 and in Section \ref{sec:16}, we consider two exterior 4-forms 
\[
 \form{\Spin{8}}, \;  \; \form{\Spin{7} \U{1}} \; \;  \in \Lambda^4 \RR^{16}  ,
\] 
canonically associated with subgroups $\Spin{8}, \, \Spin{7} \U{1} \subset \SO{16}$, and we write their explicit expressions in the 16 coordinates. We will see in  the next two statements which $4$-planes of $\RR^{16}$ are calibrated by $\form{\Spin{8}}$ and by $\form{\Spin{7} \U{1}}$. 

The 4-forms $\form{\Spin{8}}, \, \form{\Spin{7} \U{1}}$ and the respective calibrated $4$-planes can be compared with other calibrations in $\RR^{16}$, in particular with the previously mentioned Bryant-Harvey 4-form $\Phi_K$. 
It is thus appropriate to remind the main theorem in \cite[Theorem 2.27]{br-ha}, namely that, in any $\HH^n \cong \RR^{4n}$, the Bryant-Harvey $4$-form $\Phi_K$ calibrates the Cayley 4-planes that are contained in a quaternionic $2$-dimensional vector subspace $W^2_\HH \subset \HH^n$. 

Here we prove:

\begin{te}\label{Thm 1.3} The oriented 4-planes of $\RR^{16} \cong \OO^2$ calibrated by the 4-form $\form{\Spin{8}}$ are the \emph{transversal Cayley 4-planes}, i.e. the $4$-planes $P$ such that both projections $\pi(P)$, $\pi'(P)$ on the two summands in $\OO^2 =\OO \oplus \OO'$ are two dimensional and both invariant by a same complex structure $u \in S^6 \subset \Im \OO$.
\end{te}

Also, by recalling that $\OO^2$ decomposes in the union of \emph{octonionic lines} 
\[
\ell_m =\{ (x,mx), \; x \in \OO, \; m \in \OO \cup \infty\}, 
\]
meeting pairwise only at $(0,0) \in \OO^2$, we can state:

\begin{te}\label{Thm 1.4} The oriented 4-planes of $\RR^{16}$ calibrated by the 4-form $\form{\Spin{7} \U{1}}$ are the ones that are invariant under a complex structure $u \in S^6 \subset  \mathrm{Im} \OO$ and that are contained in an octonionic line $\ell_m \subset \OO \oplus \OO$, where only $m \in \RR$ and $m = \infty$ are allowed. Thus they are "Cayley 4-planes", contained in the oriented 8-planes that are the mentioned octonionic lines with $m \in \RR \cup \infty$.
\end{te}

\bigskip

\section{Preliminaries}\label{sec:prel}

The multiplication in the algebra $\OO$ of octonions can be defined from the one in quaternions $\HH$ through the Cayley-Dickson process: if $x = h_1 + h_2e, x' = h'_1 + h'_2e \in \OO$, then
\[
xx' = (h_1h'_1 - \bar h'_2h_2) + (h_2 \bar h'_1 + h'_2h_1)e,
\]
where product of quaternions is used on the right hand side and $\bar h'_1 , \bar h'_2$ are the conjugates of $h'_1 , h'_2 \in \HH$. Like for quaternions, the conjugation $\bar x = \bar h_1 - h_2 e$
in $\OO$ relates with the non-commutativity: $\overline{xx'} = \bar x' \bar x$. One has also the associator $[x, x', x''] = (xx')x'' - x(x'x'')$,
that vanishes whenever two among $x, x' , x'' \in \OO$ are equal or conjugate. 

The identification $x = h_1 + h_2e \in \OO \leftrightarrow (h_1, h_2) \in\HH^2$, used in the previous formula, is not an isomorphism of (left
or right) quaternionic vector spaces. To get an isomorphism one has instead to go through the following hypercomplex structure $(I, J, K)$ on $\RR^8 \cong \OO$. For
$x=h_1 +h_2 e \in\OO$, where $(h_1,h_2) \in\HH^2$,
define
\[
I(x) = x \cdot i, \; J(x) = x \cdot j, \;  K(x) = (x \cdot i) \cdot j
\]
or equivalently
\[
I(h_1, h_2) = (h_1i, -h_2i), \; J(h_1, h_2)) = (h_1j, -h_2j), \; K(x) = (h_1k, h_2k).
\]
This observation goes likely back to the very discovery of octonions in the mid-1800s. The alternative approach to the same isomorphism used in our Introduction does not seem however to have appeared before 1989, when R. Bryant and R. Harvey \cite{br-ha} looked at the map
\[
L: \HH^2 \rightarrow \OO, \qquad  L(h_1,h_2)=h_1 +(kh_2 \bar k)e\in\OO,
\] 
and observed it satisfies
\[
L[(h_1, h_2)i] = h_1i + (kh_2i \bar k)e, \quad L[(h_1, h_2)j] = h_1j + (kh_2j \bar k)e,\; L[(h_1, h_2)k] = h_1k + (kh_2)e.
\] 
This, in terms of $x_1 = h_1, x_2 = kh_2 \bar k$ and of the octonion $x = x_1 + x_2e$, can be read exactly as in our previous approach:
\[
L[(h_1, h_2)i] = x \cdot i, \; L[(h_1, h_2)j] = x \cdot j, \; L[(h_1, h_2)k] = (x \cdot i) \cdot j, 
\]
and as mentioned $L^*\spinform{7} = \Phi_K$.

\bigskip

\section{The quaternionic $4$-form and the Cayley calibration in $\RR^8$}\label{sec:cayley}

A possible way to produce $4$-forms canonically associated with some $G$-structures is through the notion of \emph{Clifford system}. We recall the definition, originally given in the context of isoparametric hypersurfaces, cf. \cite{fkm}.

\begin{de}A \emph{Clifford system} on a Riemannian manifold $(M,g)$ is a vector sub-bundle $E^r \subset \; \text{End}\; TM$ locally spanned by self-adjoint anti-commuting involutions $\I_1,\dots, \I_r$. Thus $I^2_\alpha=\Id,\quad\I^*_\alpha=\I_\alpha,\quad\I_\alpha\circ\I_\beta=-\I_\beta\circ\I_\alpha$, and the $\I_\alpha$ are required to be related, in the intersections of trivializing sets, by matrices of $\SO{r}$. The rank $r$ of $E$ is said to be the \emph{rank} of the Clifford system.
\end{de}

Possible ranks of irreducible Clifford systems on $\RR^N$ are classified, up to $N=32$, as follows:
\bigskip

\begin{table}[H]\label{Cliffordranks}
\caption{Rank of irreducible Clifford systems in $\RR^{N}$}
\renewcommand{\arraystretch}{2.25}
\begin{center}
{\begin{tabular}{|c|cccccccccccc|}
\hline
$ \text{dimension}\; N$&$2$&$4$&$8$&$8$&$16$&$16$&$16$&$16$&$32$&$64$&$64$& \dots \\
\hline
$\text{rank}\;  r $&$2$&$3$&$4$&$5$&$6$&$7$&$8$&$9$&$10$&$11$&$12$& \dots \\
\hline
\end{tabular}}
\end{center}
\end{table}
\bigskip

In particular, the Clifford system of rank 3 in $\RR^4$ can be defined by the classical Pauli matrices: 
\[ 
\I_1=\left(
\begin{array}{c|c}
0 & 1 \\ \hline
1 & 0
\end{array}
\right), \;
\I_2=\left(
\begin{array}{c|c}
0 & -i\\ \hline
i & 0
\end{array}
\right), \; 
\I_3=\left(
\begin{array}{c|c}
1 & 0 \\ \hline
0 & -1
\end{array}
\right) \in \U{2} \subset \SO{4},
\]
and the Clifford system of rank 5 in $\RR^8$ by the following similar \emph{(right) quaternionic Pauli matrices}:
\begin{equation}\label{qPauli}
\begin{split}
\I_1 = \left(
\begin{array}{c|c}
0 & \Id \\ \hline
\Id & 0
\end{array}
\right), \;
\I_2 = \left(
\begin{array}{c|c}
0 & -R_i \\ \hline
R_i & 0
\end{array}
\right), \;
\I_3\ =\left(
\begin{array}{c|c}
0 & -R_j \\ \hline
R_j & 0
\end{array}
\right), \;
\\
\I_4 = \left(
\begin{array}{c|c}
0 & -R_k \\ \hline
R_k & 0
\end{array}
\right), \; 
\I_5 = \left(
\begin{array}{c|c}
\Id & 0 \\ \hline
0 & -\Id
\end{array}
\right) \in \Sp{2} \subset  \SO{8},
\end{split}
\end{equation}
where as before $R_i, R_j,R_k$ denote the multiplication on the right by $i,j,k$ on $\HH^2 \cong \RR^8$.

According to Table, A, there is also a Clifford system with $r=4$ in $\RR^8$, explicitly defined by selecting \emph{e.g.}
\[
\I_1=\left(
\begin{array}{c|c}
0 & \Id \\ \hline
\Id & 0
\end{array}
\right), \;
\I_2= \left(
\begin{array}{c|c}
0 & -R_i \\ \hline
R_i & 0
\end{array}
\right), \;
\I_3 =\left(
\begin{array}{c|c}
0 & -R_j \\ \hline
R_j & 0
\end{array}
\right), \;
\I_4 =\left(
\begin{array}{c|c}
0 & -R_k\\ \hline
R_k & 0
\end{array}
\right).
\] 

Going back to rank $r=5$, from the quaternionic Pauli matrices $\I_1,\I_2,\I_3,\I_4,\I_5$, one gets the $10$ complex structures on $\RR^8$
\[
\I_{\alpha\beta} = \I_\alpha\circ\I_\beta\qquad\text{for } \qquad 1\leq \alpha<\beta \leq 5.
\]
Their K\"ahler forms $\theta_{\alpha\beta}$ give rise to a $5\times5$ skew-symmetric matrix
\[
\theta_R =(\theta_{\alpha\beta}),
\]
and one can easily see that both the following matrices of K\"ahler 2-forms 
\[
\theta_R =(\theta_{\alpha \beta}) \in \frak{so}(5) \qquad \text{and} \qquad \omega_L = \begin{pmatrix} 0&\omega_{L_i} & \omega_{L_j} \\ -\omega_{L_i} & 0 &\omega_{L_k} \\ -\omega_{L_j} &-\omega_{L_k} &0 \end{pmatrix} \in \frak{so}(3)
\]
allow to write the (left) quaternionic 4-form of $\HH^2$ as
\begin{equation}\label{thetaR}
\Omega_L= -\frac{1}{2} \sum_{\alpha < \beta} \theta^2_{\alpha \beta} =[ \omega_{L_i}^2 + \omega_{L_j}^2 + \omega_{L_k}^2 ]. 
\end{equation}

On the other hand, as mentioned in the Introduction, the subgroup $\Spin{7} \subset \SO{8}$ (generated by the right translation $R_u$, $u \in S^6 \subset \Im \OO$)
gives rise to the Cayley calibration $\Phi_{\mathrm{Spin(7)}} \in \Lambda^4$:
  \[
 \spinform{7} \; = \; -\frac{1}{6} [  \phi_i^2+  \phi_j^2+ \dots  \phi_h^2]\; = \; \frac{1}{6} \sum_{\alpha < \beta}  \zeta^2_{\alpha \beta}.
  \]
Here $ \phi_i,  \phi_j,  \phi_k,  \phi_e,  \phi_f,  \phi_g,  \phi_h.$ are the K\"ahler 2-forms associated with the complex structures $( J_1, J_2,J_3, J_4, J_5, J_6, J_7) = (R_i, R_j,R_k,R_e,R_f,R_g,R_h)$, and $\zeta=(\zeta_{\alpha \beta}) \in \frak{so}(7)$ is the matrix of the K\"ahler 2-forms of compositions $J_{\alpha \beta} = J_\alpha \circ J_\beta$.

It is worth to recall that under the action of Sp(2)Sp(1), the space of exterior 2-forms $\Lambda^2 \RR^8$ decomposes  as

 \[
\Lambda^2 = \Lambda^2_{10} \oplus \Lambda^2_{15}\oplus  \Lambda^2_{3},
\]
where lower indices denote the dimensions of irreducible components. Here
$\Lambda^2_{10} \cong \frak{sp}(2)$ is generated by the K\"ahler forms $\theta_{\alpha \beta}$ of the $\J_{\alpha \beta} (\alpha < \beta)$,
compositions of the five quaternionic Pauli matrices, and $\Lambda^2_{3} \cong \frak{sp}(1)$ is generated by the K\"ahler forms $\omega_{L_i}, \omega_{L_j}, \omega_{L_k}$.

By denoting by $\tau_2$ the second coefficient in the characteristic polynomial of the involved skew-symmetric matrices, we can rewrite formula \eqref{thetaR} of $\Omega_L$ as:

\begin{equation*}
 \Omega_L= -\frac{1}{2} \tau_2(\theta_R) =\tau_2(\omega_L)
\end{equation*}
where $\theta_R=(\theta_{\alpha \beta}) \in \mathfrak{so}(5)$, and  $\omega_L = \begin{pmatrix} 0&\omega_{L_i} & \omega_{L_j} \\ -\omega_{L_i} & 0 &\omega_{L_k} \\ -\omega_{L_j} &-\omega_{L_k} &0 \end{pmatrix} \in \mathfrak{so}(3)$.
\bigskip

Similarly, under the $\Spin{7}$ action one gets the decomposition:

\[
\Lambda^2 = \Lambda^2_{7} \oplus \Lambda^2_{21}, 
\]
where $\Lambda^2_{7} $ is generated by the K\"ahler forms $\phi_\alpha$ of the $J_\alpha = R_i,R_j, \dots , R_h$ and $\Lambda^2_{21} \cong \frak{spin}(7)$ is generated by the K\"ahler forms $\zeta_{\alpha \beta}$ of the $J_\alpha \circ J_\beta \; (\alpha < \beta)$. Thus, in the $\tau_2$ notation:

\[
\spinform{7} \; = \; -\frac{1}{6} \sum \phi_\alpha^2 \; = \; \frac{1}{6} \tau_2(\zeta), \quad \zeta=(\zeta_{\alpha \beta}) \in \mathfrak{so}(7).
\]
All the exterior 4-forms $\theta_{\alpha\beta}$, $\phi_\alpha$ and the $\zeta_{\alpha \beta}$ have been studied systematically as calibrations in the space $\RR^8$, cf. \cite{dhm}. 
\bigskip

\section{Proof of Theorems \ref{Thm 1.1} and \ref{Thm 1.2}}\label{sec:cayley2}

The matrix $\eta=(\eta_{\alpha \beta}) \in \mathfrak{so}(5)$ in the statement of Theorem \ref{Thm 1.2} is defined as follows.
Let $\I_\alpha^L \, (\alpha =1, \dots , 5)$  be the left quaternionic Pauli matrices defined as in \eqref{qPauli} but by using the left quaternionic multiplications $L_i,L_j,L_k$ by $i,j.k$. If $J^L_{\alpha \beta}= \I^L_\alpha \circ \I^L_\beta$ and if $\eta_{\alpha \beta}$ are the K\"ahler 2-forms associated to $J^L_{\alpha \beta}$, a computation shows that

\[
2 \Omega_R\; =\;  \;  \eta_{12}^2+\eta_{13}^2+\eta_{24}^2+\eta_{34}^2+\eta_{23}^2+\eta_{14}^2+\eta_{15}^2+\eta_{25}^2+\eta_{35}^2+\eta_{45}^2,
\]
and note the symmetry with the first identity in formula \eqref{thetaR}.

We express now the 2-forms $\theta_{\alpha \beta}$ and $\eta_{\alpha \beta}$ in the coordinates of $\RR^8$, using the following abridged notations. Let $\{dx_1,\dots,dx_8\}\subset\Lambda^1\RR^8$ be the standard basis of $1$-forms in $\RR^8$. Then $\shortform{\alpha\beta}$ (scriptsize) denotes $dx_\alpha\wedge dx_\beta$ and $\shortform{\alpha\beta\gamma\delta}$ denotes $dx_\alpha\wedge dx_\beta\wedge dx_\gamma\wedge dx_\delta$, and $\star$ denotes the Hodge star, so that $a+\star = a+\star a$. One gets:

\begin{equation}\label{eq:theta1}
\begin{aligned}
\theta_{12} &= -\shortform{12}+\shortform{34}+\shortform{56}-\shortform{78}\enspace, &
\theta_{13} &= -\shortform{13}-\shortform{24}+\shortform{57}+\shortform{68}\enspace, &
\theta_{14} &= -\shortform{14}+\shortform{23}+\shortform{58}-\shortform{67}\enspace,\\
\theta_{23} &= -\shortform{14}+\shortform{23}-\shortform{58}+\shortform{67}\enspace, &
\theta_{24} &= +\shortform{13}+\shortform{24}+\shortform{57}+\shortform{68}\enspace, &
\theta_{34} &= -\shortform{12}+\shortform{34}-\shortform{56}+\shortform{78}\enspace,\\
\end{aligned}
\end{equation}
and

\begin{equation}\label{eq:theta2}
\begin{aligned}
\theta_{15} &= -\shortform{15}-\shortform{26}-\shortform{37}-\shortform{48}\enspace, &
\theta_{25} &= -\shortform{16}+\shortform{25}+\shortform{38}-\shortform{47}\enspace, \\
\theta_{35} &= -\shortform{17}-\shortform{28}+\shortform{35}+\shortform{46}\enspace, &
\theta_{45} &= -\shortform{18}+\shortform{27}-\shortform{36}+\shortform{45}\enspace,
\end{aligned}
\end{equation}
so that, if $\theta=(\theta_{\alpha\beta})$ 

\begin{equation}\label{Theta}
\begin{split}
\tau_2(\theta)=\theta^2_{12} + \theta^2_{13} + \dots + \theta^2_{45} = \hspace{4cm}\\ 
= -12\shortform{1234}-4\shortform{1256}-4\shortform{1357}+4\shortform{1368}-4\shortform{1278}-4\shortform{1467}-4\shortform{1458}+\star = -2\Omega_L.
\end{split}
\end{equation}
Next:

\begin{equation}\label{eq:eta1}
\begin{aligned}
\eta_{12} &= -\shortform{12}-\shortform{34}+\shortform{56}+\shortform{78}\enspace, &
\eta_{13} &= -\shortform{13}+\shortform{24}+\shortform{57}-\shortform{68}\enspace, &
\eta_{14} &= -\shortform{14}-\shortform{23}+\shortform{58}+\shortform{67}\enspace,\\
\eta_{23} &= +\shortform{14}+\shortform{23}+\shortform{58}+\shortform{67}\enspace, &
\eta_{24} &= -\shortform{13}+\shortform{24}-\shortform{57}+\shortform{68}\enspace, &
\eta_{34} &= +\shortform{12}+\shortform{34}+\shortform{56}+\shortform{78}\enspace,\\
\eta_{15} &= -\shortform{15}-\shortform{26}-\shortform{37}-\shortform{48}\enspace, &
\eta_{25} &= -\shortform{16}+\shortform{25}-\shortform{38}+\shortform{47}\enspace, \\
\eta_{35} &= -\shortform{17}+\shortform{28}+\shortform{35}-\shortform{46}\enspace, &
\eta_{45} &= -\shortform{18}-\shortform{27}+\shortform{36}+\shortform{45}\enspace,
\end{aligned}
\end{equation}
and, if $\eta=(\eta_{\alpha\beta})$, 

\begin{equation}\label{Eta}
\begin{split}
\tau_2(\eta)=\eta^2_{12} + \eta^2_{13} + \dots + \eta^2_{45} = \hspace{4cm}\\ 
= 12\shortform{1234}-4\shortform{1256}-4\shortform{1357}-4\shortform{1368}+4\shortform{1278}+4\shortform{1467}-4\shortform{1458}+\star = -2\Omega_R.
\end{split}
\end{equation}
Similarly:

\begin{equation}\label{phi7}
\begin{aligned}
\phi_{i} = -\shortform{12}+\shortform{34}+\shortform{56}-\shortform{78}\enspace,\qquad
\phi_{j} &= -\shortform{13}-\shortform{24}+\shortform{57}+\shortform{68}\enspace,\qquad &
\phi_{k} &= -\shortform{14}+\shortform{23}+\shortform{58}-\shortform{67}\enspace,\qquad \\
\phi_{e} &= -\shortform{15}-\shortform{26}-\shortform{37}-\shortform{48}\enspace,\qquad &
\phi_{f} &= -\shortform{16}+\shortform{25}-\shortform{38}+\shortform{47}\enspace,\qquad \\
\phi_{g} &= -\shortform{17}+\shortform{28}+\shortform{35}-\shortform{46}\enspace,\qquad &
\phi_{h} &= -\shortform{18}-\shortform{27}+\shortform{36}+\shortform{45}\enspace,
\end{aligned}
\end{equation}
and one easily deduce also formulas for the $\zeta_{\alpha\beta}$ (cf. \cite{pp1}, \cite{Boy}). Then, by \eqref{spin7}:

\begin{equation}\label{Spin7}
\spinform{7} = \shortform{1234}+\shortform{1256}+\shortform{1357}+\shortform{1368}-\shortform{1278}-\shortform{1467}+\shortform{1458}+\star\enspace.
\end{equation}
By computing the squares of the 2-forms in \eqref{phi7} \eqref{eq:eta1} and comparing with Formulas \eqref{Eta}, \eqref{Theta}, \eqref{Spin7}, the identities listed in Theorems \ref{Thm 1.1} and \ref{Thm 1.2} are recognized.
\bigskip

\section{Even Clifford structures in dimension $8$}

We recall first the following notion, proposed in 2001 by A. Moroianu and U. Semmelmann, \cite{ms}.

\begin{de} Let $(M,g)$ be a Riemannian manifold. An \emph{even Clifford structure} is the choice of an oriented Euclidean vector bundle $E^r$ of rank $r \geq 2$ over $M$, together with a bundle morphism $\varphi$ from the even Clifford algebra bundle
\begin{equation*}
\varphi: \mathrm{Cl}^{even} E^r \rightarrow \mathrm{End} TM \quad \text{such that} \quad \Lambda^{2} E^r \hookrightarrow \mathrm{End^-} TM.
\end{equation*}
$r$ is called the  {\emph rank}  of the even Clifford structure.
\end{de}

The even Clifford structure $E$ is said to be \emph{parallel} if there exists a metric connection $\nabla^E$ on $E$ such that $\varphi$ is connection preserving, i.e.
\[
\varphi(\nabla^E_X \sigma) = \nabla^g_X \varphi (\sigma),
\]
for every tangent vector $X \in TM$ and section $\sigma$ of $\mathrm{Cl}^{even} E$, where $\nabla^g$ is the Levi Civita connection. 

Rank $2,3,4$ parallel even Clifford structures are equivalent to complex K\"ahler, quaternion K\"ahler, product of two quaternion K\"ahler. Besides them, higher rank parallel non-flat even Clifford structures in dimension $8$ are listed in the following Table, cf. \cite{ms}

\begin{table}[H]
\caption{Parallel non-flat even Clifford structures of rank $\geq 5$ in $M^8$}
\begin{center}
\renewcommand{\arraystretch}{1.65}
\begin{tabular}{|c||c|}
\hline
$r$ &$M$\\
\hline\hline
5&\scriptsize{quaternion K\"ahler}\\
\hline
6&\scriptsize{K\"ahler} \\
\hline
7&\scriptsize{$\mathrm{Spin}(7)$} \scriptsize{holonomy}\\
\hline
8&\scriptsize{Riemannian}\\
\hline
\end{tabular}
\end{center}
\end{table}

A class of examples of even Clifford structures are those coming from Clifford systems as defined in Section \ref{sec:cayley}. Namely, if 
the vector sub-bundle $E^r \subset \; \text{End} \; TM$, locally spanned by self-adjoint anticommuting involutions $\I_1,\dots, \I_r$, defines the Clifford system, then one easily recognizes that through the compositions $J_{\alpha \beta} = \I_\alpha \circ \I_\beta$, the Clifford morphism $\varphi: \mathrm{Cl}^{even}(E^r) \rightarrow \mathrm{End} \; TM$ is well defined.

An example is given by the first row of the former Table, where the quaternion K\"ahler structure is constructed via the local $J_{\alpha \beta} = \I_\alpha \circ \I_\beta$ defined as in Section \ref{sec:cayley}, by using on the model space $\HH^2$ the quaternionic Pauli matrices. The remaining three rows of the former Table correspond to \emph{essential} even Clifford structures, i. e. to even Clifford structures that cannot be defined throw a Clifford system, cf. \cite{ppv} for a discussion on this notion.

The following Table gives a description of the Clifford bundle generators and of the canonically associated 4-form for each of the four even Clifford structures on $\RR^8$.

\begin{table}[H]
\caption{Generators and associated 4-forms}
\begin{center}
\renewcommand{\arraystretch}{1.95}
\begin{tabular}{|c||c|c|c|}
\hline
r &$M$ &Clifford bundle generators& associated 4-form \\
\hline\hline
5&qK&$\I_1,\I_2,\I_3,\I_4,\I_5$& $\Omega_L = -\frac{1}{2}\tau_2 (\theta_{\alpha \beta}), \quad \Omega_R = -\frac{1}{2}\tau_2 (\eta_{\alpha \beta})$\\
\hline\hline
6&K\"ahler &$J_1,J_2,J_3,J_4,J_5,J_6$& $\Phi_{\Spin{6}}= \tau_2(\zeta_{\alpha \beta})=-5\omega^2$\\
\hline
7&$\Spin{7}$ hol& $J_1,J_2,J_3,J_4,J_5,J_6,J_7$&$\spinform{7}=-\frac{1}{6} \sum \phi_\alpha^2 =\frac{1}{6} \tau_2(\zeta_{\alpha \beta})$\\
\hline\hline
8&Riemannian&$I,J_1,J_2,J_3,J_4,J_5,J_6 ,J_7$&$\Phi_{\SO{8}}= \tau_2(\psi_{\alpha \beta})= 0$ \\
\hline
\end{tabular}

\end{center}
\end{table}

Here $\I_1,\I_2,\I_3,\I_4,\I_5$ are the (left or right) quaternionic Pauli matrices, $ (\theta_{\alpha \beta})$, $(\eta_{\alpha \beta})$ are like in Section \ref{sec:cayley2}. Notations $(J_1,J_2,J_3,J_4,J_5,J_6 ,J_7)=(R_i,R_j,R_k,R_e,R_f,R_g,R_h)$ are also used, $\phi_\alpha$ is the K\"ahler form of $J_\alpha$ and $\zeta_{\alpha \beta}$ is the K\"ahler form of $J_\alpha \circ J_\beta$. Finally, $(\psi_{\alpha \beta}) \in \mathfrak{so}(8)$, with entries $\pm \phi_\alpha$ in the first line and column  and with entries $\zeta_{\alpha \beta} \in \mathfrak{so}(7)$.

It is of course desirable to give examples of Riemannian manifolds $(M^8,g)$ supporting both a $\Sp{2}\cdot \Sp{1}$ and a $\Spin{7}$ structure. Rarely the metric $g$ can be the same for both structures, but this is possible of course for parallelizable $(M^8,g)$. On this respect, homogeneous $(M^8,g)$ with an invariant $\Spin{7}$ structure have been recently classified \cite{acfr}, by making use of the following topological condition for compact oriented spin $M^8$:
\[
p_1^2(M) -4 p_2(M) + 8\chi(M) =0.
\]
Some of the obtained examples are parallelizable, \emph{e. g.} diffeomorphic to $S^7 \times S^1$ and $S^5 \times S^3$. On the latter, $S^5 \times S^3$, using two natural parallelizations, one can define two $\Spin{7}$ structures, both of general type (in the 1986 M. Fernandez $\Spin{7}$ framework), and the hyperhermitian structure associated with one of them corresponds to a family of Calabi-Eckmann \cite{pp4}.

To get examples of $8$-dimensional manifolds that admit both a locally conformally hyperk\"ahler metric $g$ and a locally conformal parallel $\Spin{7}$ metric, that is either the same $g$ as before, or a different metric $g'$, a good point to start with is the class of compact 3-Sasakian 7-dimensional manifolds $(\mathcal S^7,g)$. Many examples of such $(\mathcal S^7,g)$ and with arbitrary second Betti numbers have been given by Ch. Boyer -K. Galicki \emph{et al}, cf. \cite{bgmr}. In particular, recall that given the 3-Sasakian $(\mathcal S^7,g)$ one gets a locally conformally hyperk\"ahler metric $g$ on the product $\mathcal S^7 \times S^1$ \cite{op}. This can also be expressed by saying that the 3-Sasakian metric $g$  has the property of being nearly parallel $\Gtwo$, and in particular with $3$ linearly independent Killing spinors, cf. \cite[pages 536-538]{bg}. Moreover the differentiable manifold $\mathcal S^7$ admits, besides the $3$-Sasakian metric $g$, another metric $g'$ that is also nearly parallel $\Gtwo$ but \emph{proper, i. e.} with only one non zero Killing spinor. This allows to extend the metrics $g$ and $g'$ to the product with $S^1$ and to get both the properties of locally conformally hyperk\"ahler and locally conformally parallel $\Spin{7}$ on $(M^8, g)=(\mathcal S^7 \times S^1, g)$ and of locally conformally parallel $\Spin{7}$ on $(M^8, g' )=(\mathcal S^7 \times S^1, g')$, cf. also \cite{ipp}.

Further examples of $8$-dimensional  differentiable manifolds admitting both a $\Sp{2}\cdot \Sp{1}$-structure with respect to a metric $g$ and a $\mathrm{Spin}(7)$-structure with respect to a metric $g'$ include the Wolf spaces $\HH P^2$ and $\Gtwo / \SO{4}$, cf. \cite{acfr}.  Finally, the non singular sextic 
$Y =\{[z_0, \dots ,z_5] \in \CC P^5, \; z_0^6+ \dots +z_5^6 =0\}$ is also an example, where a metric $g$ giving an almost quaternionic structure is insured by a result in \cite{cv}, and a metric $g'$ with holonomy $\SU{4} \subset \Spin{7}$ by Calabi-Yau theorem, cf \cite[p. 139]{j}.
\bigskip

\section{Dimension $16$}\label{sec:16}

A Clifford system with $r=9$ in $\OO^2 \cong \RR^{16}$ is given by the following \emph{octonionic Pauli matrices}: 

\[
\I_1 \hspace{-0.1cm}= \hspace{-0.1cm}\left(
\begin{array}{c|c}
0 & \Id \\ \hline
\Id & 0
\end{array}
\right), \;
\I_2 \hspace{-0.1cm}= \hspace{-0.1cm}\left(
\begin{array}{c|c}
0 & -R_i \\ \hline
R_i & 0
\end{array}
\right), \; 
\I_3 \hspace{-0.1cm}= \hspace{-0.1cm}\left(
\begin{array}{c|c}
0 & -R_j \\ \hline
R_j & 0
\end{array}
\right),  
\] 
\[
\I_4 \hspace{-0.1cm}= \hspace{-0.1cm}\left(
\begin{array}{c|c}
0 & -R_k \\ \hline
R_k & 0
\end{array}
\right), \; 
\I_5  \hspace{-0.1cm}=\hspace{-0.1cm}\left(
\begin{array}{c|c}
0 & -R_e \\ \hline
R_e & 0
\end{array}
\right),  \; 
\I_6 \hspace{-0.1cm}=\hspace{-0.1cm}\left(
\begin{array}{c|c}
0 & -R_f\\ \hline
R_f & 0
\end{array}
\right), 
\] 
\[ 
\I_7 \hspace{-0.1cm}=\hspace{-0.1cm}\left(
\begin{array}{c|c}
0 & -R_g \\ \hline
R_g & 0
\end{array}
\right), \;
\I_8 \hspace{-0.1cm}=\hspace{-0.1cm}\left(
\begin{array}{c|c}
0 & -R_h \\ \hline
R_h & 0
\end{array}
\right), \;
\I_9 \hspace{-0.1cm}= \hspace{-0.1cm}\left(
\begin{array}{c|c}
\Id & 0 \\ \hline
0 & -\Id
\end{array}
\right) \in \SO{16},
\]
and of course now $R_i,R_j, \dots , R_h$ denote the multiplication on the right by the unit octonions $i,j,\dots, h$ on $\OO^2 \cong \RR^{16}$.

Looking back at Table A, we see that in $\RR^{16}$ there are also irreducible Clifford systems with $r=8,7,6$. According to \cite{ppv}, convenient choices are the following:
\[
r=8: \; 
  \I_1, \dots, \I_8, \qquad 
r=7: \;
 \I_2 ,\dots,
\I_8,
\qquad
r=6: \;
 \I_1 ,\I_2, \I_3,
\I_4, 
\I_5, 
\I_9.
\]

It is now worth to remind the following parallel situations in complex, quaternionic and octonionic geometry. The groups $\U{2}\subset\SO{4}$, $\Sp{2}\cdot\Sp{1}\subset\SO{8}$, $\Spin{9} \subset \SO{16}$ are the stabilizers of the vector subspaces
\[
E^{3} \subset \mathrm{End^+}  (\RR^4), \qquad E^{5} \subset \mathrm{End^+}  (\RR^8), \qquad E^{9} \subset \mathrm{End^+} (\RR^{16})
\]
spanned respectively by the Pauli, quaternionic Pauli, octonionic Pauli matrices.

Moreover, $ \U{2}$, $\Sp{2}\cdot\Sp{1}$, $\Spin{9}$ are symmetry groups of the Hopf fibrations respectively:
\[
S^{3}\overset{S^1}\longrightarrow S^2\cong\CP{1}, \qquad S^{7}\overset{S^3}\longrightarrow S^4\cong\HP{1}, \qquad  S^{15}\overset{S^7}\longrightarrow S^8\cong\OP{1}.
\]

Finally, $\U{2}$, $\Sp{2}\cdot\Sp{1}$, $\Spin{9}$ are stabilizers in $\Lambda^2 \CC^2 , \Lambda^4 \HH^2, \Lambda^8 \OO^2$ of the following \emph{canonically associated forms}, cf. \cite{BerCCP}:
\begin{equation*}
\scriptsize{\form{\U{2}}\hspace{-0.1cm}= \hspace{-0.1cm}\int_{\CP{1}}  p_\ell^*\nu_\ell d\ell \in \Lambda^2, \quad 
\form{\Sp{2}\cdot\Sp{1}} \hspace{-0.1cm} = \hspace{-0.1cm} \int_{\HP{1}}p_\ell^*\nu_\ell\,d\ell \in \Lambda^4, \quad 
\form{\Spin{9}} \hspace{-0.1cm}= \hspace{-0.1cm} \int_{\OP{1}}p_\ell^*\nu_\ell\,d\ell \in \Lambda^8},
\end{equation*}
where $\nu_\ell$ is the volume form on the line $\ell \ug\{(x,mx)\}$ or $\ell \ug\{(0,y)\}$ in $\CC^2$ or $\HH^2$ or $\OO^2$,
$$p_\ell:\; \CC^2 \cong \RR^4 \;  \text{or} \; \HH^2 \cong \RR^8 \; \text{or} \; \OO^2 \cong \RR^{16} \longrightarrow \ell$$ 
is the projection on the line $\ell)$, and note that the integral formula is based on the volume of distinguished planes.
In the three cases one gets in this way the K\"ahler 2-form of $\CC^2$, the quaternion K\"ahler 4-form of $\HH^2$ 
and the canonical 8-form of $\OO^2$.
\bigskip

\section{Rank $8$, $7$ and $6$ Clifford systems on $\RR^{16}$}\label{r=8,7}

Look now closer at the nine octonionic Pauli matrices, that define a rank $9$ Clifford system in $\RR^{16}$, and at the choices among them that give rise to ranks $r=8,7,6$ (cf. previous Section). 
The compositions $\I_{\alpha \beta}=\I_\alpha \circ \I_\beta, \; \alpha<\beta$,
for all choices $r=6,7,8,9$ are bases of the Lie algebras
\[
\mathfrak{spin}(6) \subset \mathfrak{spin}(7) \subset \mathfrak{spin}(8) \subset \mathfrak{spin}(9) \subset \mathfrak{so}(16).
\]
Like in the previous Sections, we can write the matrices of K\"ahler forms $\psi_{\alpha \beta}$ associated to $\I_{\alpha \beta}$, and we use for them the following notations:
\[
\psi^A =(\psi_{\alpha \beta}) \in \mathfrak{so}(6), \; \; \psi^B=(\psi_{\alpha \beta}) \in \mathfrak{so}(7), \; \; \psi^C =(\psi_{\alpha \beta}) \in \mathfrak{so}(8), \; \; \psi^D =(\psi_{\alpha \beta}) \in \mathfrak{so}(9).
\]
The second coefficients $\tau_2$ of their characteristic polynomial give rise to the following invariant $4$-forms
\[
 \; \tau_2(\psi^A), \; \; \tau_2(\psi^B), \; \; \tau_2(\psi^C), \; \; \tau_2(\psi^D) \; \;  \in \; \;  \Lambda^4 \RR^{16}\; 
\]
that can be written in (the differentials of) the coordinates of $\RR^{16} = \OO \oplus \OO$:
\[
\shortform{1,2,3,4,5,6,7,8;1',2',3',4',5',6',7',8'}.
\]
We recall in particular that in the $\Spin{9}$ situation, the following identity holds:
\[
\tau_2(\psi^D) =0,
\]
and this gives evidence to the next coefficient $\tau_4(\psi^D) \in \Lambda^8$, proportional to the $8$-form $\form{\Spin{9}}$, as studied in \cite{pp1}.
\bigskip

\section{The $4$-forms $\spinform{8}$ and $\form{\Spin{7} \U{1}}$}

Look now only at $\I_1, \dots , \I_8$ and at the matrix 
$$\psi^C=(\psi_{\alpha \beta}) \in \mathfrak{so}(8)$$
of K\"ahler forms associated to $ \I_{\alpha \beta}=\I_\alpha \circ \I_\beta.$ By using coordinate 1-forms $$\shortform{1,2,3,4,5,6,7,8;1',2',3',4',5',6',7',8'},$$ an explicit computations based on the explicit formulas for the $\psi_{\alpha \beta}$ in \cite{pp1} yields:

\[
\boxed{
\begin{array}{rl}
&\spinform{8}=\frac{1}{4}\tau_2(\psi^C)=\frac{1}{4}\displaystyle\sum^8_{1=\alpha<\beta}\psi^2_{\alpha\beta}=\boldshortform{1234}+\boldshortform{1256}-\boldshortform{1278}+\boldshortform{1357}+\boldshortform{1368}+\boldshortform{1458}\\
&-\boldshortform{1467}-\boldshortform{2358}+\boldshortform{2367}+\boldshortform{2457}+\boldshortform{2468}-\boldshortform{3456}+\boldshortform{3478}+\boldshortform{5678} +\displaystyle\sum^8_{1=a<b} \shortform{aba'b'} \; \in \; \Lambda^4 \RR^{16}.
\end{array}
}
\]
\bigskip

\noindent Here boldface notations have the following meaning:
\[
\boldshortform{abcd} =\shortform{abc'd'}-\shortform{ab'cd'}+\shortform{ab'c'd}+\shortform{a'bcd'}-\shortform{a'bc'd}+\shortform{a'b'cd}.
\]
\bigskip

By excluding now the K\"ahler forms involving $\I_1$ and $\I_9$, one gets the matrix $\psi^B =(\psi_{\alpha \beta}) \in \mathfrak{so}(7)$. Now similar computations lead to:

\[
\boxed{
\begin{array}{rl}
&\qquad \qquad \qquad \form{\Spin{7} \U{1}}=\tau_2(\psi^B)=  6[\shortform{1234}+\shortform{1256}\\
&-\shortform{1278}+\shortform{1357}+\shortform{1368}+\shortform{1458}-\shortform{1467}-\shortform{2358}+\shortform{2367}+\shortform{2457}+\shortform{2468}-\shortform{3456}+\shortform{3478}+\shortform{5678}]\\
&+6[\shortform{1'2'3'4'}+\shortform{1'2'5'6'}-\shortform{1'2'7'8'}+\shortform{1'3'5'7'}+\shortform{1'3'6'8'}+\shortform{1'4'5'8'}-\shortform{1'4'6'7'}-\shortform{2'3'5'8'}\\
&\qquad +\shortform{2'3'6'7'}+\shortform{2'4'5'7'}+\shortform{2'4'6'8'}-\shortform{3'4'5'6'}+\shortform{3'4'7'8'}+\shortform{5'6'7'8}]+6\displaystyle\sum^8_{1=a<b}\shortform{aba'b'}\\
&+2[\boldshortform{1234}+\boldshortform{1256}-\boldshortform{1278}+\boldshortform{1357}+\boldshortform{1368}+\boldshortform{1458}-\boldshortform{1467}-\boldshortform{2358}\\
&\qquad \qquad+\boldshortform{2367}+\boldshortform{2457}+\boldshortform{2468}-\boldshortform{3456}+\boldshortform{3478}+\boldshortform{5678}] \; \in \; \mathrm{\Lambda^4 \RR^{16}},
\end{array}
}
\]
\bigskip

\noindent where boldface notations have the same meaning as before. The presence of the factor $\U{1}$ in the group $\Spin{7}\U{1}$ is here due to a computation showing that matrices in $\SO{16}$ commuting with the seven involutions $\I_2, \dots \I_8$ are a $\U{1}$ subgroup, well identified in \cite[Chapter 6, p. 44]{Boy} (cf. also proof of Theorem \ref{Thm 1.4}  below).

We take this opportunity to remark that the expression of $\form{\Spin{7} \U{1}}$ written in the paper \cite{ppv} contains some errors in the coefficients as well as some missing terms, and has to be corrected by the present one. Note also that $\form{\Spin{7} \U{1}}$ restricts, on any of the two summands of $\RR^{16}=\RR^8 \oplus \RR^8$, and up to a factor $6$, to the usual Cayley calibration of \cite{ha-la}.
\bigskip


We are now ready for the proofs of Theorems \ref{Thm 1.3} and \ref{Thm 1.4}. The following notion has already implicitly introduced in the statement of Theorem \ref{Thm 1.3}.
\bigskip

\begin{de} Let $P$ be a $4$-plane in the real vector space $\RR^{16} \cong \OO^2 = \OO \oplus \OO'$, and let $\pi: \OO^2 \rightarrow \OO$ and $\pi': \OO^2 \rightarrow \OO'$ be the orthogonal projections to $\OO$ and $\OO'$. $P$ is said to be a \emph{transversal Cayley $4$-plane} if both its projections $\pi(P)$, $\pi'(P)$ are $2$-dimensional and invariant under a same complex structure $u \in S^6 \subset \Im \OO$.
\end{de}
\bigskip

\noindent \emph{Proof of Theorem \ref{Thm 1.3}}. Recall that $\Spin{8}$ can be characterized as the subgroup of the following matrices $A\in \SO{16}$:
\[ 
A= \left(
\begin{array}{c|c}
a_+&0\\
\hline
0&a_- 
\end{array}
\right)
\]
where $a_+,a_- \in \SO{8}$ are \emph{triality companions}, i. e. and for any $v \in \OO$ there exists a $w \in \OO$ such that $R_w =a_+R_v a_-^t$ (cf. \cite[p. 278-279]{ha}). It follows that $\Spin{8}$ contains the diagonal $\Spin{7}_\Delta$ (characterized by choices $a_+=a_-)$ and acts transitively on transversal $4$-planes of $\RR^{16}$. On the other hand the 4-form $\form{\Spin{8}}$ is invariant under the action of $\Spin{8}$. Thus, since $\form{\Spin{8}}$ takes value $1$ on the $4$-plane spanned by the coordinates $\shortform{121'2'}$, $\form{\Spin{8}}$ takes value $1$ on any tranversal Cayley 4-plane in $\RR^{16}$. 

Next, let $Q$ be any $4$-plane of $\RR^{16}$. By looking at the expression of $\form{\Spin{8}}$, we see that the only possibilities for having non zero value on $Q$ are that $\pi(Q)$ and $\pi'(Q)$ are $2$-dimensional. For such $4$-planes $Q$ we can use the following canonical form with respect to the complex structure $i \in S^6$:
\[
Q = \big[ e_1 \wedge (R_i e_1 \cos \theta + e_2 \sin \theta) \big] \oplus \big[ e'_1 \wedge (R_i e'_1 \cos \theta' + e'_2 \sin \theta') \big],
\]
where the pairs $e_1, e_2$ and $e'_1, e'_2$ are both orthonormal and respectively in $\OO$ and in $\OO'$, and with angles limited by $0 \leq \theta \leq \frac{\pi}{2}$ and $\theta \leq \theta' \leq \pi- \theta$. The above canonical form for $Q$ is a small variation of the canonical forms that are used in a proof of the classical Wirtinger's inequaliy (cf. \cite[p.6]{l}) and in characterizations of Cayley $4$-planes in $\RR^8$ in the Harvey-Lawson foundational paper (cf. \cite[p. 121]{ha-la}). Its proof follows the steps of proof of the mentioned canonical form, as explained in details in \cite{l}. From this canonical form we see that $\form{\Spin{8}}(Q) \leq 1$ for any $4$-plane $Q$, and that the equality holds only if $\theta = \theta' =0$, i. e. for transversal Cayley $4$-planes. 
\bigskip

\noindent \emph{Proof of Theorem \ref{Thm 1.4}}. The leading terms in the expression of $\form{\Spin{7} \U{1}}$ are those with coefficient $6$, thus terms involving only coordinates among $\shortform{12345678}$, or only coordinates among $\shortform{1'2'3'4'5'6'7'8'}$, or terms $\shortform{aba'b'}$. Look first at the first and second types of terms. We already mentioned that the restriction of $\form{\Spin{7} \U{1}}$ to any of the summands in $\OO^2 =\OO \oplus \OO'$ is the usual Cayley calibration in $\RR^8$, whose calibrated $4$-planes are the Cayley planes. Thus, for the first two types of terms, we get as calibrated $4$-planes just the Cayley $4$-planes that are contained in the octonionic lines with slope $m=0$ and $m=\infty$. In the remaining case of terms $\shortform{aba'b'}$ one gets as calibrated $4$-planes the transversal Cayley $4$-planes that are contained in the octonionic line $\ell_1$ (leading coefficient $m=1$). Now $\Spin{7}$ acts on the individual octonionic lines $\ell_0, \ell_1, \ell_\infty$, and the only possibility to move planes out of them is through the factor $\U{1}$. In fact, the discussion in \cite[Chapter 6, p. 44]{Boy} shows that the factor $\U{1}$ in the group $\Spin{7}\U{1}$ moves the octonionic lines through the circle, contained in the space $S^8$ of the octonionic lines, passing through the three points $m=0,1,\infty$. This corresponds to admitting any real coefficient: $m \in \RR \cup \infty$ as slope of the octonionic lines that are admitted to contain the calibrated $4$-planes.
\bigskip

\begin{re} Following the recent work \cite{Kor} by J. Kotrbat\'y, one can use octonionic 1-forms, according to the following formal definitions:

\begin{equation*}
\begin{split}
dx = d\alpha+ id\beta+ j d\gamma +kd\delta +e d\epsilon +f d\zeta +g d\eta +h d\theta, \\
\overline{dx} = d\alpha- id\beta- j d\gamma -kd\delta -e d\epsilon -f d\zeta -g d\eta -h d\theta, \\
dx' = d\alpha'+ id\beta'+ j d\gamma' +kd\delta' +e d\epsilon' +f d\zeta' +g d\eta' +h d\theta', \\ 
\overline{dx'} = d\alpha'- id\beta'- j d\gamma' -kd\delta' -e d\epsilon' -f d\zeta' -g d\eta' -h d\theta',
\end{split}
\end{equation*}
referring to pairs of octonions $(x,x') \in \OO \oplus \OO = \RR^{16}$. Then, in the same spirit proposed in \cite{Kor}, a straightforward computation yields the following formula, much simpler way to write the $\Spin{8}$ canonical 4-form of $\RR^{16}$:

\[
\; \spinform{8} = \frac{1}{4}(\overline{dx} \wedge dx') \wedge (\overline{dx'} \wedge dx).\; 
\]
\bigskip

\noindent Similarly, one gets that the $\Spin{7}\U{1}$ canonical  $4$-form of $\RR^{16}$ can be written in octonionic 1-forms as:

\begin{equation*}
\begin{split}
\; \form{\Spin{7}\U{1}} = \frac{1}{4} \big[ (\overline{dx} \wedge dx)^2 + (\overline{dx'} \wedge dx')^2 \big] 
\hspace{4.5cm} 
\\
\; 
\hspace{2.5cm} -\frac{1}{2} \big[ (\overline{dx} \wedge dx')^2+(\overline{dx'} \wedge dx)^2 \big] -\big[(\overline{dx} \wedge dx')\wedge (\overline{dx'} \wedge dx)\big]. \;
\end{split}
\end{equation*}
Details of both computations are in \cite{Boy}.

\end{re}

\end{document}